\documentclass{amsart}
\usepackage[lite]{amsrefs}
\usepackage{amssymb,amsmath,amsthm,amscd,latexsym}
\usepackage{amssymb}
\usepackage{amsfonts}
\usepackage[mathscr]{eucal}
\usepackage{color}
\usepackage{graphicx}
\usepackage{float}
\usepackage[all,cmtip]{xy}
\usepackage{tikz}
\usepackage{placeins}

\newcommand{\FF}{\mathbb{F}}

\newcommand{\ds}{\displaystyle}

\begingroup
\newtheorem{thrm}{Theorem\hskip 1mm}[section]
\newtheorem{crry}[thrm]{Corollary\hskip 1mm}
\newtheorem{prop}[thrm]{Proposition\hskip 1mm}

\newtheorem{lemm}[thrm]{Lemma\hskip 1mm}

\endgroup

\newenvironment{pf}{\noindent\textbf{Proof.}}{\hspace*{\fill}$\blacksquare$\\[6pt]}

\begin{document}

\title{New Self-dual Codes from $2 \times 2$ block circulant matrices, Group Rings and Neighbours of Neighbours}

\author{J. Gildea, A. Kaya, R. Taylor and A. Tylyshchak}

\address{Department of Mathematics\\
Faculty of Science and Engineering\\
University of Chester\\
England}
\email{j.gildea@chester.ac.uk}

\address{Department of Mathematics Education, Sampoerna University\\
12780, Jakarta, Indonesia}
\email{nabidin@gmail.com}

\address{Department of Mathematics\\
Faculty of Science and Engineering\\
University of Chester\\
England}
\email{rhian.taylor@chester.ac.uk}

\address{Department of Algebra\\
Uzhgorod National University\\
Uzhgorod, Ukraine}
\email{alxtlk@bigmir.net}


\begin{abstract}

In this paper, we construct self-dual codes from a construction that involves $2 \times 2$ block circulant matrices, group rings and a reverse circulant matrix. We provide
conditions whereby this construction can yield self-dual codes. We construct
self-dual codes of various lengths over $\FF_2$, $\FF_2+u\FF_2$ and $\FF_4+u\FF_4$. Using extensions, neighbours
and neighbours of neighbours, we construct $32$ new self-dual codes of length $68$.

\end{abstract}

\maketitle

\section{Introduction}
 Self-dual codes (a class of linear block codes) is a vibrant area of Mathematics which was first introduced in the early 1970's (\cite{BE,AG,MMS,MST}).
 The double circulant construction (introduced in \cite{CPW,MK}) is one of the most extensively used techniques to construct self-dual codes. It involves
 considering generator matrices of the form $(I|A)$ where $A$ is a circulant matrix. In \cite{GL}, the authors consider constructing self-dual codes
 from generator matrices of the form $(I|A)$ where $A$ is a block circulant matrix. In \cite{BLM}, certain well known self-dual codes were constructed
 from certain group rings. In recent years group rings have been used to construct self-dual codes in \cite{Us,Us2}. In this article we construct
 self-dual codes by considering generator matrices that combine $2 \times 2$ block circulant construction, group rings and reverse circulant matrices. In particular
 we construct self-dual codes from generator matrices of the form:\\

 \begin{center}
\begin{tikzpicture}[scale=0.8]
\node at (0,0) {$A$}; \node at (1.5,0) {$B+C$};
\node at (0,-0.75) {$B+C$}; \node at (1.5,-0.75) {$A$};
\draw (-1,0.2)--(-1,-1);
\draw (2.4,0.2)--(2.5,0.2)--(2.5,-1)--(2.4,-1);
\node at (-2,-0.375) {$I$};
\draw (-2.9,0.2)--(-3,0.2)--(-3,-1)--(-2.9,-1);
\end{tikzpicture}
\end{center}

where $A$ and $B$ are matrices that arise from a group ring construction and $C$ is a reverse circulant matrix. For the remainder of this section, we will introduce many
important concepts and results required for further sections. In section 2, we describe the construction itself. We present the structure of the generator matrix and discuss associated theory in order to put some restrictions on unknowns. These restrictions aim to maximise the practicality of the construction method by reducing the search field. Following the theory, we look at the numerical results from certain groups of order $4$, $8$ and $17$. We then apply extensions and consider neighbours of codes as methods of finding new codes.\\

Throughout this article, we assume that $R$ is a finite Frobenius ring of characteristic $2$.  A code over a finite commutative ring  $R$ is defined as any subset $C$ of $R^n$. An element of $C$ is called a codeword. If a code satisfies $C=C^\perp$ then the code $C$ is said to be self-dual, alternatively if $C \subseteq C^\perp$ then the code is said to be self-orthogonal. The Hamming weight enumerator of a code is defined as:
\[
W_C(x,y) = \sum_{\mathbf{c} \in C} x^{n-wt(\mathbf{c})} y^{wt(\mathbf{c})}.
\]

\noindent For binary codes, a self-dual code where all weights are congruent to $0 \pmod{4}$ is said to be Type II and the code is said to be Type I otherwise. If a code satisfies $W_C(x,y) = W_{C^\perp}(x,y)$ then the code is said to be formally self-dual. The bounds on the minimum distances, $d(n)$ for Type I and  Type II codes respectively, are (\cite{ER})
\[ d(n) \leq 4 \lfloor \frac{n}{24} \rfloor +4 \]
and
\[ d(n) \leq \begin{cases} 4 \lfloor \frac{n}{24} \rfloor +4 ~\text{if}~ n \not \equiv 22~(\text{mod}\;24) \\
4 \lfloor \frac{n}{24} \rfloor +6 ~\text{if}~ n \equiv 22~(\text{mod}\;24) \end{cases}
 \]
If these bounds are met for self-dual codes, they are called extremal. Although, all of the theoretical results are based around finite Frobenius rings of characteristic $2$,
all of the numerical results are based on the rings $\FF_2$, $\FF_2+u\FF_2$ and $\FF_4+u\FF_4$.

The first commutative ring that we consider is $\FF_2+u\FF_2 := \FF_2 [X]/(X^2)$, where $u$ satisfies $u^2=0$. The elements of the ring may be written as $0,1,u$ and $1+u$, where $1$ and $1+u$ are the units of $\FF_2+u\FF_2$. We also consider $\FF_4+u\FF_4$; the commutative binary ring of size $16$. $\FF_4+u\FF_4$ can be viewed as an extension of $\FF_2+u\FF_2$. Therefore, we can express any element of $\FF_4 +u\FF_4$ in the form $\omega a+(1+\omega)b$, where $a, b\in \FF_2 + u\FF_2$. These rings are generalised in \cite{Rkpaper} and \cite{Rkpaper2}. The most effective way of displaying these results, is to use the hexadecimal system. This is achieved by use of the ordered basis $\{u\omega, \omega, u, 1\}$.
\begin{eqnarray*}
0 &\leftrightarrow &0000,\ 1\leftrightarrow 0001,\ 2\leftrightarrow 0010,\
3\leftrightarrow 0011, \\
4 &\leftrightarrow &0100,\ 5\leftrightarrow 0101,\ 6\leftrightarrow 0110,\
7\leftrightarrow 0111, \\
8 &\leftrightarrow &1000,\ 9\leftrightarrow 1001,\ A\leftrightarrow 1010,\
B\leftrightarrow 1011, \\
C &\leftrightarrow &1100,\ D\leftrightarrow 1101,\ E\leftrightarrow 1110,\
F\leftrightarrow 1111.
\end{eqnarray*}

The following Gray Maps were introduced in \cite{gaborit,ling} and \cite%
{dougherty2};%
\begin{eqnarray*}
\psi _{\mathbb{F}_{4}} &:&a\omega +b\overline{\omega }\mapsto \left(
a,b\right) \text{, \ }a,b\in \mathbb{F}_{2}^{n} \\
\varphi _{\mathbb{F}_{2}+u\mathbb{F}_{2}} &:&a+bu\mapsto \left( b,a+b\right)
\text{, \ }a,b\in \mathbb{F}_{2}^{n} \\
\psi _{\mathbb{F}_{4}+u\mathbb{F}_{4}} &:&a\omega +b\overline{\omega }%
\mapsto \left( a,b\right) \text{, \ }a,b\in \left( \mathbb{F}_{2}+u\mathbb{F}%
_{2}\right) ^{n} \\
\varphi _{\mathbb{F}_{4}+u\mathbb{F}_{4}} &:&a+bu\mapsto \left( b,a+b\right)
\text{, \ }a,b\in \mathbb{F}_{4}^{n}
\end{eqnarray*}%
These Gray maps preserve orthogonality in the respective alphabets, for the
details we refer to \cite{QRFpvFp,ling}. The binary codes $\varphi _{\mathbb{F}_{2}+u%
\mathbb{F}_{2}}\circ \psi _{\mathbb{F}_{4}+u\mathbb{F}_{4}}\left( C \right) $
and $\psi _{\mathbb{F}_{4}}\circ \varphi _{\mathbb{F}_{4}+u\mathbb{F}%
_{4}}\left( C\right) $ are equivalent to each other.

\begin{prop}
$($\cite{ling}$)$ Let $C$ be a code over $\mathbb{F}_{4}+u\mathbb{F}_{4}$.
If $C$ is self-orthogonal, so are $\psi _{\mathbb{F}_{4}+u\mathbb{F}%
_{4}}\left( C\right) $ and $\varphi _{\mathbb{F}_{4}+u\mathbb{F}_{4}}\left(
C\right) $. $C$ is a Type I (resp. Type II) code over $\mathbb{F}_{4}+u%
\mathbb{F}_{4}$ if and only if $\varphi _{\mathbb{F}_{4}+u\mathbb{F}%
_{4}}\left( C\right) $ is a Type I (resp. Type II) $\mathbb{F}_{4}$-code, if
and only if $\psi _{\mathbb{F}_{4}+u\mathbb{F}_{4}}\left( C\right) $ is a
Type I (resp. Type II) $\mathbb{F}_{2}+u\mathbb{F}_{2}$-code. Furthermore,
the minimum Lee weight of $C$ is the same as the minimum Lee weight of $\psi
_{\mathbb{F}_{4}+u\mathbb{F}_{4}}\left( C\right) $ and $\varphi _{\mathbb{F}%
_{4}+u\mathbb{F}_{4}}\left( C\right) $.
\end{prop}

\begin{crry}
Suppose that $C$ is a self-dual code over $\mathbb{F}_{4}+u\mathbb{F}_{4}$
of length $n$ and minimum Lee distance $d$. Then $\varphi _{\mathbb{F}_{2}+u%
\mathbb{F}_{2}}\circ \psi _{\mathbb{F}_{4}+u\mathbb{F}_{4}}\left( C\right) $
is a binary $\left[ 4n,2n,d\right] $ self-dual code. Moreover, $C$ and $%
\varphi _{\mathbb{F}_{2}+u\mathbb{F}_{2}}\circ \psi _{\mathbb{F}_{4}+u%
\mathbb{F}_{4}}\left( C\right) $ have the same weight enumerator. If $C$ is
Type I (Type II), then so is $\varphi _{\mathbb{F}_{2}+u\mathbb{F}_{2}}\circ
\psi _{\mathbb{F}_{4}+u\mathbb{F}_{4}}\left( C\right) $.
\end{crry}

\begin{thrm}
\label{ext} (\cite{Kim}) Let $C$ be a self-dual code of length $n$
over a commutative\ Frobenius ring with identity $R$ and $G=(r_i)$ be a $k \times n$ generator matrix for $C$, where $r_i
$ is the i-th row of $G $, $1\leq i \leq k.$ Let $c$ be a unit in $R$ such
that $c^2=-1$ and $X$ be a vector in $S^n$ with $\langle X,X \rangle=-1.$
Let $y_i=\langle r_i, X \rangle $ for $1 \leq i \leq k.$ The following matrix

\begin{equation*}
\begin{bmatrix}
\begin{tabular}{cc|c}
$1$ & $0$ & $X$ \\ \hline
$y_1$ & $cy_1$ & $r_1$ \\
$\vdots$ & $\vdots$ & $\vdots$ \\
$y_k$ & $cy_k$ & $r_k$%
\end{tabular}%
\end{bmatrix}%
,
\end{equation*}
generates a self-dual code $D$ over $R$ of length $n+2.$
\end{thrm}

Two self-dual binary codes of length $2n$ are said to be neighbours
of each other if their intersection has dimension $n-1$. Let $x\in {\mathbb{F}}%
_{2}^{2n}-\mathcal{C}$ then $\mathcal{D}=\left\langle \left\langle
x\right\rangle ^{\bot }\cap \mathcal{C},x\right\rangle $ is a neighbour of $%
\mathcal{C}$. \\

 The definitions surrounding group rings are as follows: Let $G$ be a finite group of order $n$, then the group ring $RG$ consists of $\sum_{i=1}^n \alpha_i g_i$, $\alpha_i \in R$, $g_i \in G$. Addition in the group ring is done by coordinate addition:
\begin{equation} \sum_{i=1}^n \alpha_i g_i + \sum_{i=1}^n \beta_i g_i =
\sum_{i=1}^n (\alpha_i + \beta_i ) g_i.\end{equation}  The product of two elements in a group ring  is defined as:
\begin{equation} \left( \sum_{i=1}^n \alpha_i g_i \right) \left( \sum_{j=1}^n \beta_j g_j\right)  = \sum_{i,j} \alpha_i \beta_j g_i g_j. \end{equation}  It follows, that the coefficient of $g_i$ in the product is $ \sum_{g_i g_j = g_k } \alpha_i \beta_j .$ Note that, $e_G$ denotes the identity element of any group $G$.

The following construction of a matrix was first given for codes over fields by Hurley, \cite{Hurley1}, and extended to rings in \cite{Us}. Let $R$ be a finite commutative Frobenius ring and let $G = \{ g_1,g_2,\dots,g_n \}$ be the elements of a group of order $n$ in a given listing. Let $v = \sum_{i=1}^n \alpha_{g_i}   \in RG.$  Define the matrix $\sigma(v) \in M_n(R)$ to be $\sigma(v)=(\alpha_{g_i^{-1} g_j})$ where $i,j \in \{1,2,\cdots,n\}$. \\
\noindent The main group discussed in this work is the cyclic group. A circulant $n\times n$ matrix is denoted $cir(\alpha_1, \alpha_2, \cdots, \alpha_n)$, where each row vector is rotated one element to the left relative to the preceding row vector \cite{davis}. Additionally, a reverse circulant $n\times n$ matrix is denoted $rcir(\alpha_1, \alpha_2, \cdots, \alpha_n)$, where each row vector is rotated one element to the right relative to the preceding row vector. The notation $CIR(A_1, A_2, \cdots, A_m)$ denotes the block circulant matrix where the first row of block matrices are $A_1,\ldots,A_n$.  If $v=\ds{ \sum_{i=0}^{n-1}} \alpha_{i}x^i\in RC_{n}$, then $\sigma(v)=\mbox{circ}(\alpha_0,\alpha_1,\dots,\alpha_{n-1})$
where $\alpha_i \in R$. We will now look at the structure of the matrix $\sigma(v)$ where $v$ is an element of $C_{2p}$. \\

\noindent Let $C_{2p}= \langle x~|~x^{2p}=1 \rangle$ and
\[ v=\sum_{i=0}^{p-1} \sum_{j=0}^{1} \alpha_{i+pj+1} x^{2i+j} \in RC_{2p} \]
\noindent then,
\[ \sigma(v)=\begin{pmatrix} A_1 & A_2 \\
A_2' & A_1 \end{pmatrix} \]
\noindent where $A_j=cir(\alpha_{(j-1)p+1}, \alpha_{(j-1)p+2}, \dots, \alpha_{jp})$ and $A'_j=cir(\alpha_{jp}, \alpha_{(j-1)p+1}, \dots, \alpha_{jp-1})$. \\

Recall the canonical involution $* : RG \rightarrow RG$ on a group ring $RG$ is given by $v^*=\sum_{g}\alpha_g g^{-1},$ for $v=\sum_g \alpha_g g \in RG.$ If $v$ satisfies $vv^*=1,$ then we say that $v$ is a unitary unit in $RG.$ We also note that $\sigma(v^*)=\sigma(v)^T$.

\section{Construction}

Consider the following matrix $M(\sigma)$, where $v_1$ and $v_2$ are distinct group ring elements from the same group ring $RG$ where $R$ is a finite Frobenius commutative ring of characteristic $2$ and $G$ is a finite group of order $n$. $\sigma(v)$ is a matrix generated from a group ring element and $A$ denotes a reverse circulant matrix.

 \begin{center}
\begin{tikzpicture}[scale=0.8]
\node at (0,0) {$\sigma(v_1)$}; \node at (2,0) {$\sigma(v_2)+A$};
\node at (0,-0.75) {$\sigma(v_2)+A$}; \node at (2,-0.75) {$\sigma(v_1)$};
\draw (-1,0.2)--(-1,-1);
\draw (2.9,0.2)--(3,0.2)--(3,-1)--(2.9,-1);
\node at (-2,-0.375) {$I_{2n}$};
\draw (-2.9,0.2)--(-3,0.2)--(-3,-1)--(-2.9,-1);
\node at (-4.3,-0.375) {$M(\sigma)=$};
\end{tikzpicture}
\end{center}

\noindent Let $C_{\sigma}$ be the code generated by the matrix $M(\sigma)$. Clearly, $C_{\sigma}$ has length $4n$. We will now establish conditions when $C_{\sigma}$ generates a self-dual code.
We will also create a link between unitary units in $RG$ and the above construction yielding self-dual codes.

\begin{lemm}\label{ll1}
Let $R$ be a finite commutative Frobenius ring of characteristic $2$ and let $B$ and $C$ be $n \times n$ matrices over $R$. Then, the matrix
 \begin{center}
\begin{tikzpicture}[scale=0.8]
\node at (0,0) {$B$}; \node at (1,0) {$C$};
\node at (0,-0.75) {$C$}; \node at (1,-0.75) {$B$};
\draw (-1,0.2)--(-1,-1);
\draw (1.9,0.2)--(2,0.2)--(2,-1)--(1.9,-1);
\node at (-2,-0.375) {$I_{2n}$};
\draw (-2.9,0.2)--(-3,0.2)--(-3,-1)--(-2.9,-1);
\node at (-4.3,-0.375) {$M=$};
\end{tikzpicture}
\end{center}
generates a self-dual code iff $(B+C)(B+ C)^T=I_n$ and $BC^T=C B^T$.
\end{lemm}

\begin{pf} Clearly, the code generated by $M$ has free rank $2n$, as the left-hand side of the matrix $M$ is the $2n \times 2n$ identity matrix. The code generated by $M$
is self-dual iff the code generated by $M$ is self-orthogonal. Now,

\[
MM^T
=I_{2n}+
\left(\begin{array}{cc}
B&C\\
C&B\\
\end{array}\right)
\left(\begin{array}{cc}
B^T&C^T\\
C^T&B^T\\
\end{array}\right)=\left(\begin{array}{cc}
I_n+BB^T+ C C^T&BC^T+ C B^T\\
CB^T+BC^T&I_n+ C C^T+BB^T\\
\end{array}\right)
\]
and $MM^T=0$ iff $I_n+BB^T+ C C^T=0$ and $BC^T+C B^T=0$. Adding these equations, we obtain
\[ I_n+BB^T+ C C^T+BC^T+ C B^T=0 \Longleftrightarrow (B+C)(B+ C)^T=I_n.\]

\end{pf}

\begin{thrm}\label{th:sd}
Let $R$ be a finite commutative Frobenius ring of characteristic $2$ and let $G$ be a finite group of order $n$. Then, $C_{\sigma}$ generates a self-dual code of length $4n$ iff
$(\sigma(v_1+v_2)+A)(\sigma((v_1+v_2)^*)+A)=I_n$ and $\sigma(v_1)(\sigma( (v_1+v_2)^*)+A)=(\sigma(v_1+v_2)+A)\sigma(v_1^*)$.
\end{thrm}
\begin{pf}
By the previous result, $C_{\sigma}$ generates a self-dual code iff
\[(\sigma(v_1)+\sigma(v_2)+A)(\sigma(v_1)+\sigma(v_2)+A)^T=I_n \;\text{and} \;\sigma(v_1)(\sigma(v_2)+A)^T=(\sigma(v_2)+A)\sigma(v_1)^T.\]

\noindent Now, $\sigma(v_1)+\sigma(v_2)+A=\sigma(v_1+v_2)+A$ and
\[
\begin{split}
(\sigma(v_1)+\sigma(v_2)+A)^T &= \sigma(v_1)^T+\sigma(v_2)^T+A^T\\
&=\sigma(v_1^*)+\sigma(v_2^*)+A\\
&=\sigma(v_1^*+v_2^*)+A\\
&=\sigma((v_1+v_2)^*)+A.
\end{split}
\]
\noindent Clearly, $\sigma(v_1)(\sigma(v_2)+A)^T=(\sigma(v_2)+A)\sigma(v_1)^T$ is equivalent to
\[\sigma(v_1)\sigma(v_1)^T+\sigma(v_1)(\sigma(v_2)+A)^T=\sigma(v_1)\sigma(v_1)^T+(\sigma(v_2)+A)\sigma(v_1)^T.\]
\noindent Considering the left-and right-hand sides separately, we obtain:
\[
\begin{split}
 \sigma(v_1)\sigma(v_1)^T+\sigma(v_1)(\sigma(v_2)+A)^T&=\sigma(v_1)\sigma(v_1^*)+ \sigma(v_1)(\sigma(v_2)^T+A^T)\\
&=\sigma(v_1)\sigma(v_1^*)+ \sigma(v_1)\sigma(v_2^*)+\sigma(v_1)A\\
&=\sigma(v_1)(\sigma(v_1^*)+ \sigma(v_2^*)+A)\\
&= \sigma(v_1)(\sigma(v_1^*+v_2^*)+A)\\
&= \sigma(v_1)(\sigma((v_1+v_2)^*)+A).
\end{split}
\]
\noindent and
\[
\begin{split}
(\sigma(v_2)+A)\sigma(v_1)^T+\sigma(v_1)\sigma(v_1)^T &=\sigma(v_1)\sigma(v_1^*)+ (\sigma(v_2)+A)\sigma(v_1^*)\\
&=\sigma(v_1)\sigma(v_1^*)+ \sigma(v_2)\sigma(v_1^*)+A\sigma(v_1^*)\\
&= (\sigma(v_1)+\sigma(v_2)+A)\sigma(v_1^*)\\
&= (\sigma(v_1+v_2)+A)\sigma(v_1^*).
\end{split}
\]
\end{pf}

\begin{lemm}\label{ll3}
Let $R$ be a finite commutative Frobenius ring of characteristic $2$, $A$ be an $n\times n$ reverse
circulant over $R$ and $V$ be an $n\times n$ circulant matrix over $R$.
Then,
\begin{equation}\label{ee4}
AV^T+VA^T=0.
\end{equation}
\end{lemm}

\begin{pf} Let $V=\text{circ}(v_1,v_n,v_{n-1},\ldots,v_3,v_2)$. Clearly, $V=v_1I_n+v_2P+v_3P^2+\cdots+v_nP^{n-1}$ where
$P=\text{circ}(0,0,\ldots,0,1)$ and $A=\text{rcirc}(a_1,a_2,\ldots,a_{n-1},a_n)$. Now,
\[
\begin{split}
V^T&=v_1I_n^T+v_2P^T+v_3(P^2)^T+\cdots+v_n(P^{n-1})^T\\
&=v_1I_n+v_2P^T+v_3(P^T)^2+\cdots+v_n(P^T)^{n-1}.
\end{split}
\]

\noindent As $A=A^T$, it remains to show that $AP^T+PA=0$. Finally,

\[
PA=\text{circ}(0,0,\ldots,0,1) \cdot \text{rcirc}(a_1,a_2,\ldots,a_{n-1},a_n)=\text{rcirc}(a_n,a_1,\ldots,a_{n-1})
\]
\noindent and
\[
AP^T=\text{rcirc}(a_1,a_2,\ldots,a_{n-1},a_n) \cdot \text{circ}(0,1,\ldots,0,0)=\text{rcirc}(a_n,a_1,\ldots,a_{n-1}).
\]

\end{pf}

\begin{lemm}\label{l:sym}
Let $R$ be a commutative ring and let $G=\{g_1=e,\ldots,g_n\}$ be a finite group of order $n>1$.
The $\sigma(v)$ is symmetric for any $v\in RG$ if and only if $G$ is abelian group of exponent $2$.
\end{lemm}

\begin{pf}
Clearly, $\sigma(v)$ is symmetric for any $v\in RG$ if and only if
$\alpha_{g_i^{-1}g_j}=\alpha_{g_j^{-1}g_i}$ $(i,j=1,\ldots,n)$ for any $v=\sum_{g\in G}\alpha_g g\in RG.$
Furthermore, we have $g_i^{-1}g_j=g_j^{-1}g_i$ $(i=1,\ldots,n)$ $(i,j=1,\ldots,n)$ or
$xy=y^{-1}x^{-1}$ for any $x,y\in G$.
Note that for an abelian group of exponent $2$, $yxy=x^{-1}$ or $xyxy=e$  or $(xy)^2=e$ for any $x,y\in G$. Therefore, we have that $g^2=e$ for any $g\in G$; thus, $G$ has exponent $2$.

It is interesting to note that any group of exponent $2$ is abelian because $xyxy=e$ and $xxyy=ee=e$
since $x$ and $y$ are commutative for any $x,y\in G$.
\end{pf}

\begin{lemm}
Let $R$ be a commutative ring.
An $n\times n$-matrix $X$ satisfies $XA=AX^T$ for any $n\times n$ reverse
circulant matrix $A$ over $R$ if and only if $X$ is a circulant matrix.
\end{lemm}

\begin{pf}
This proof follows from lemma \ref{ll3}. Let $X$ be
an $n\times n$-matrix which satisfies $XA=AX^T$.
Then
\[XA=A^TX^T\]
and
\[XA=(XA)^T\]
for any $n\times n$ reverse circulant matrix $A$ over $R$.
This implies that $XA$ is symmetric. Let $D=\begin{pmatrix}0&\ldots&0&1\\0&\ldots&1&0\\\vdots&\ddots&\vdots&\vdots\\1&\ldots&0&0\end{pmatrix}=\text{rcirc}(0,0,\ldots,0,1)$,
 $X=(x_{i,j})$.
Clearly, we have $D^2=I_n$ and $XDDA$ is symmetric for any $n\times n$ reverse circulant matrix $A$ over $R$. Therefore, $(x_{{i},{n-j}})DA$ is symmetric.

So we have $(x_{{i},{n-j}})B$ is symmetric for any $n\times n$
circulant matrix $B$ over $R$.
This is equivalent to the fact that $(x_{{i},{n-j}})P^k$ is symmetric for any $k\in\{1,\ldots,n\}$ and
$n\times n$ matrix $P=\text{circ}(0,0,\ldots,0,1)$.
Thus, $(x_{{i},{(k-j) \mbox{\footnotesize mod} n+1}})$ is symmetric for any $k\in\{1,\ldots,n\}$.
We have
\[
x_{{i},{(k-j) \mbox{\footnotesize\, mod\,} n+1}}=x_{{j},{(k-i) \mbox{\footnotesize\, mod\,} n+1}}\ i,j,k\in\{1,\ldots,n\}
\]

It is easy to see that $j'=(k-j) \mbox{\footnotesize\, mod\,} (n+1)$ equivalent to  $j=(k-j') \mbox{\footnotesize\, mod\,} (n+1)$ where $i,j,j',k\in\{1,\ldots,n\}$.
So
\[
x_{{i},{j'}}=x_{{(k-j') \mbox{\footnotesize\, mod\,} n+1},{(k-i) \mbox{\footnotesize \,mod\,} n+1}}\ i,j',k\in\{1,\ldots,n\}
\]
Thus $((k-j')\mbox{\footnotesize\, mod\,} (n+1)) -((k-i) \mbox{\footnotesize\, mod\,} (n+1)) \equiv i - j \pmod{ n}.$
Therefore, we have that $x_{{i},{j'}}$ is constant if $(i-j)\, \mbox{mod}\, n$ is fixed. Thus, $X$ is circulant.
\end{pf}

\begin{lemm}
Let $R$ be a finite commutative Frobenius ring of characteristic $2$ and let $G$ be a finite abelian group of order $n$ of exponent $2$.
Then, $C_{\sigma}$ generates a self-dual code of length $4n$ if
$\sigma(v_1), \sigma(v_2)$ are circulant matrices,
$\sigma((v_1+v_2)^2)+A^2=I_n$.
\end{lemm}

\begin{pf} We note that $A\sigma(v_1^*)=\sigma(v_1)A$, $A\sigma(v_2^*)=\sigma(v_2)A$ by lemma \ref{ll3}.
By lemma \ref{l:sym} for any $v\in RG$ $\sigma(v)$ is symmetric, so $\sigma(v^*)=\sigma(v)^T=\sigma(v).$
We also know by theorem \ref{th:sd} that $C_{\sigma}$ generates a self-dual code iff
\[(\sigma(v_1)+\sigma(v_2)+A)(\sigma(v_1)+\sigma(v_2)+A)^T=I_n \;\text{and} \;\sigma(v_1)(\sigma(v_2)+A)^T=(\sigma(v_2)+A)\sigma(v_1)^T.\]
\noindent Now,
\[
\begin{split}
(\sigma(v_1)+\sigma(v_2)+A)(\sigma(v_1)+\sigma(v_2)+A)^T&=(\sigma(v_1+v_2)+A)(\sigma((v_1+v_2)^*)+A)\\
&=\sigma(v_1+v_2)\sigma((v_1+v_2)^*)+[\sigma(v_1+v_2)A+A\sigma((v_1+v_2)^*)]+A^2\\
&=\sigma((v_1+v_2)(v_1+v_2)^*)+A^2=\sigma((v_1+v_2)^2)+A^2=I_n.
\end{split}
\]
\noindent and
\[
\begin{split}
\sigma(v_1)(\sigma(v_2)+A)^T+(\sigma(v_2)+A)\sigma(v_1)^T&=\sigma(v_1)\sigma(v_2^*)+[\sigma(v_1)A+A\sigma(v_1^*)]+\sigma(v_2)\sigma(v_1^*)\\
&=\sigma(v_1v_2)+\sigma(v_2v_1)\\
&=\sigma(v_1v_2)+\sigma(v_1v_2)=0.
\end{split}
\]
\end{pf}

\begin{lemm}
Let $R$ be a finite commutative Frobenius ring of characteristic $2$ and let $G$ be a finite cyclic group of order $n$. Then, $C_{\sigma}$ generates a self-dual code of length $4n$ iff
$\sigma((v_1+v_2)(v_1+v_2)^*)+A^2=I_n$ and $v_1v_2^*=v_2v_1^*$.
\end{lemm}

\begin{pf} We note that $A\sigma(v^*)=\sigma(v)A$ for all $v \in RG$ by the previous result. We also know that $C_{\sigma}$ generates a self-dual code iff
\[(\sigma(v_1)+\sigma(v_2)+A)(\sigma(v_1)+\sigma(v_2)+A)^T=I_n \;\text{and} \;\sigma(v_1)(\sigma(v_2)+A)^T=(\sigma(v_2)+A)\sigma(v_1)^T.\]
\noindent Now,
\[
\begin{split}
(\sigma(v_1)+\sigma(v_2)+A)(\sigma(v_1)+\sigma(v_2)+A)^T&=(\sigma(v_1+v_2)+A)(\sigma((v_1+v_2)^*)+A)\\
&=\sigma(v_1+v_2)\sigma((v_1+v_2)^*)+[\sigma(v_1+v_2)A+A\sigma((v_1+v_2)^*)]+A^2\\
&=\sigma((v_1+v_2)(v_1+v_2)^*)+A^2=I_n
\end{split}
\]
\noindent and
\[
\begin{split}
\sigma(v_1)(\sigma(v_2)+A)^T+(\sigma(v_2)+A)\sigma(v_1)^T&=\sigma(v_1)\sigma(v_2^*)+[\sigma(v_1)A+A\sigma(v_1^*)]+\sigma(v_2)\sigma(v_1^*)\\
&=\sigma(v_1v_2^*)+\sigma(v_2v_1^*)\\
&=\sigma(v_1v_2^*+v_2v_1^*).
\end{split}
\]
\noindent Finally, $\sigma(v_1v_2^*+v_2v_1^*)=0$ iff $v_1v_2^*=v_2v_1^*$.
\end{pf}

\begin{lemm}
Let $R$ is a finite commutative Frobenius ring of characteristic $2$ and let $G$ be a finite abelian group of order $n$. Let $C_{\sigma}$ be self-dual. If $A=0$, then $v_1+v_2$ is unitary.
\end{lemm}
\begin{pf}
If $C_{\sigma}$ is self-dual and $A=0$, then $\sigma((v_1+v_2)(v_1+v_2)^*)=I_n$ and $(v_1+v_2)(v_1+v_2)^*=1$.
\end{pf}

This concludes the theoretical part of this paper. We will now show the numerical results.

\section{Numerical Results}

\noindent In this section, we construct 32 new self-dual codes of length $68$. We begin with the construction of self-dual codes of length $64$ from
groups of order $4$ and $8$. Using Theorem \ref{ext}, we construct new self-dual codes of length $68$. Next, we construct codes of length $68$
from groups of order $17$. We then find new self-dual codes of length $68$ by finding neighbours of these codes and neighbours of these neighbours. Magma (\cite{magma}) was used to construct all of the codes throughout this section.\\

The possible weight enumerators for a self-dual Type I $\left[ 64,32,12\right]$-code are given in \cite{conway,binary} as:
\begin{eqnarray*}
W_{64,1} &=&1+\left( 1312+16\beta \right) y^{12}+\left( 22016-64\beta
\right) y^{14}+\cdots ,14\leq \beta \leq 284, \\
W_{64,2} &=&1+\left( 1312+16\beta \right) y^{12}+\left( 23040-64\beta
\right) y^{14}+\cdots ,0\leq \beta \leq 277.
\end{eqnarray*}%

Extremal singly even self-dual codes with weight enumerators $W_{64,1}$ are known (\cite{anev,YA,HANKEL}):

\[
\beta \in \left\{\begin{array}{c}
14, 16, 18, 19,20, 22, 24, 25, 26, 28, 29, 30, 32, 34, \\
35, 36, 38, 39, 44, 46, 49, 53, 54, 58, 59, 60, 64, 74
\end{array}\right\}
\]

and extremal singly even self-dual codes with weight enumerator $W_{64,2}$ are known for:

\[
\beta \in
\left\{\begin{array}{c}
0, . . . ,40, 41, 42, 44, 45, 46, 47, 48, 49, 50, 51, 52, 54, 55, 56, 57,  \\
58, 60, 62, 64, 69, 72, 80, 88, 96, 104, 108, 112, 114, 118, 120, 184
\end{array}\right\} \setminus \{31, 39\}.
\]

The weight enumerator of a self-dual $[68,34,12]_{2}$ code is in one of the
following forms:
\[
\begin{split}
W_{68,1}&=1+(442+4\beta )y^{12}+(10864-8\beta )y^{14}+\dots ,\\
W_{68,2}&=1+(442+4\beta )y^{12}+(14960-8\beta -256\gamma )y^{14}+\dots \ ,
\end{split}
\]
where $\beta $ and $\gamma $ are parameters and $0\leq \gamma \leq 9.$\\

 The existence of codes in $W_{68,1}$ are known for (\cite{Us2}) $%
\beta =$104, 105, 112, 115, 117, 119, 120, 122, 123, 125,\ldots , 284, 287, 289,291,
294, 301, 302, 308, 313, 315, 322, 324, 328,\ldots , 336, 338, 339, 345,
347, 350, 355, 379 and 401.\\

The first examples of codes with a $\gamma =7$ in $W_{68,2}$ are constructed in
\cite{YIL}.  Together with these, the existence of the codes in $W_{68,2}$
is known for the following parameters (see \cite{YIL,HANKEL}):

$%
\begin{array}{l}
\gamma =0,\ \beta \in \{2m|m=0,7,11,14,17,21,\dots
,99,102,105,110,119,136,165\};\ \text{or} \\
\beta \in \{2m+1|m=3,5,8,10,15,16,17,20,\dots ,82,87,93,94,101,104,110,115\};
\\
\gamma =1,\ \beta \in \{2m|m=19,22,\dots ,99\};\ \text{or}\ \beta \in
\{2m+1|m=24,\dots ,85\}; \\
\gamma =2,\ \beta \in \{2m|m=29,\dots ,100,103,104\};\ \text{or} \beta \in\{2m+1|m=32,\dots ,81,84,85,86\}; \\
\gamma  =6\text{ with }\beta \in \left\{ 2m|m=69,77,78,79,81,88\right\}  \\
\gamma  =7\text{ with }\beta \in \left\{ 7m|m=14,\ldots,39,42\right\} .
\end{array}%
$\newline

Firstly, we construct self-dual codes of length $64$ from $C_4$ (over $\FF_4+u\FF_4$), $C_{4,2}$ (over $\FF_2+u\FF_2$) and $C_{8}$ (over $\FF_2+u\FF_2$). We then construct three self-dual
codes of length $68$ (Table 4) by applying theorem \ref{ext} to the codes constructed in Tables 1,2 and 3. We replace $1+u \in \FF_2+u\FF_2$ with $3$ to save space.

\FloatBarrier
\begin{table}[H]\label{T_1}
\caption{Self-dual code over $\FF_4+u\FF_4$ of length $32$ from $C_{4}$ and $C_{4}$.}
\begin{center}
\scalebox{0.8}{
\begin{tabular}{|c|ccccc|}
\hline
$A_{i}$& $v \in C_{4}$ & $v \in C_{4}$&$r_A$ & $|Aut(A_i)|$ & $\beta$  \\ \hline \hline
 $1$ &$(8966)$ & $(0000)$&$(A617)$ & $2^4$ & $0$\\ \hline
\end{tabular}}%
\end{center}
\end{table}
\FloatBarrier

\begin{table}[h!]\label{T_2}
\caption{Self-dual code over $R_1$ of length $64$ from $C_{8}$ and $C_{8}$.}
\begin{center}
\scalebox{0.8}{
\begin{tabular}{|c|ccccc|}
\hline
$B_{i}$& $v \in C_{8}$ & $v \in C_{8}$&$r_A$ & $|Aut(B_i)|$ & $\beta$ \\ \hline \hline
$1$ & $(uuu10311)$ & $(uu011uu0)$&$(u0300013)$ & $2^3$ & $0$  \\ \hline
\end{tabular}}%
\end{center}
\end{table}

\begin{table}[h!]\label{T_3}
\caption{Self-dual code over $R_1$ of length $64$ from $C_{42}$ and $C_{42}$.}
\begin{center}
\scalebox{0.8}{
\begin{tabular}{|c|ccccc|}
\hline
$C_{i}$& $v \in C_{42}$ & $v \in C_{42}$&$r_A$ & $|Aut(C_i)|$ & $\beta$ \\ \hline \hline
$1$ & $(uu01u0u1)$ & $(u0u11u31)$&$(u3u3u3u0)$ & $2^{4}$ & $48$ \\ \hline
\end{tabular}}%
\end{center}
\end{table}

\begin{table}[h!]\label{T_4}
  \begin{center}
    \caption{Self-dual code of length $68$ from extensions of $C_1$, $C_2$ and $C_3$.}
  \scalebox{0.8}{ \begin{tabular}{|c|cccccc|} \hline
$D_i$&  Code & $c$  & $X$ & $\gamma$  & $\beta$ &$|Aut(E_i)|$ \\ \hline
$1$ &  $A_1$ &$1$ & $(0133010303011u1001333u01031uuu1u)$ & $4$ & $113$ & $2$\\ \hline
$2$ &  $B_1$ &$u+1$ & $(013011030003013301111030uuu13u10)$ & $\textbf{2}$ & $\textbf{61}$ & $2$\\ \hline
$3$ &  $C_1$ &$u+1$ & $(0u10303u110333001103u00130103303)$ & $\textbf{1}$ & $\textbf{179}$ & $2$\\ \hline
    \end{tabular}}
  \end{center}
\end{table}

\noindent We now construct two self-dual codes of length $68$ using $C_{17}$ (Table 5). We let $v_2=0 \in RC_{17}$. We note that in this case, the construction is equivalent
to the usual four circulant construction.

\begin{table}[h!]
\caption{Self-dual codes over $\FF_2$ of length $68$ $(W_{682})$ from $C_{17}$ and $C_{17}$.}
\begin{center}
\scalebox{0.8}{
\begin{tabular}{|c|cccccc|}
\hline
$E_{i}$& $v_1 \in C_{17}$ & $v_2 \in C_{17}$ & $r_A$ & $|Aut(D_i)|$ & $\gamma$ & $\beta$  \\ \hline \hline
$1$& (00000000000011011) &$(00000000000000000)$ & $(00100110010110111)$ & $2^2 \cdot 17$ & $0$ & $238$  \\ \hline \hline
$2$& (00000000110001111) & $(00000000000000000)$ & $(00100100101010101)$ & $2^2 \cdot 17$ & $0$ & $272$  \\ \hline
\end{tabular}}%
\end{center}
\end{table}

\noindent We now construct neighbours of these codes and neighbours of these neighbours.

\begin{table}[h]\caption{New codes of length 68 from neighbours of $E_1$ and $E_2$}
\begin{center}\scalebox{0.8}{
\begin{tabular}{ccccccc}
\hline
$F_{i}$ & $E_{i}$ & $(x_{35},x_{36},...,x_{68})$ & $|Aut(F_i)   |$ & $\gamma$ & $\beta$ & Type \\ \hline
$1$ & $2$ & $(0111011100100011000001001000100110)$ & $2$ & $\textbf{0}$ & $\textbf{208}$ & $W_{68,2}$  \\ \hline
$2$ & $2$ & $(1110000011111000011000011110011000)$ & $1$ & $\textbf{0}$ & $\textbf{214}$ & $W_{68,2}$  \\ \hline
$3$ & $2$ & $(0001000100001110111100001010011010)$ & $2$ & $\textbf{1}$ & $\textbf{191}$ & $W_{68,2}$  \\ \hline
$4$ & $2$ & $(0010111111111110001111001010111001)$ & $2$ & $\textbf{1}$ & $\textbf{202}$ & $W_{68,2}$  \\ \hline
$5$ & $1$ & $(1001101111101110011000101000010110)$ & $1$ & $\textbf{1}$ & $\textbf{210}$ & $W_{68,2}$  \\ \hline
$6$ & $2$ & $(0101001000111001100011110011000101)$ & $1$ & $\textbf{1}$ & $\textbf{211}$ & $W_{68,2}$  \\ \hline
$7$ & $2$ & $(0010101101010100111100000001010001)$ & $1$ & $\textbf{1}$ & $\textbf{229}$ & $W_{68,2}$  \\ \hline
$8$ & $2$ & $(1111111111111111111011101111111111)$ & $1$ & ${}$ &  $\textbf{317}$ & $W_{68,1}$  \\ \hline
\end{tabular}}
\end{center}
\end{table}

\begin{table}[h!]\caption{New codes of length 68 from neighbours of $F_7$ and $F_8$}
\begin{center}\scalebox{0.8}{
\begin{tabular}{ccccccc}
\hline
$G_{i}$ & $F_{i}$  & $(x_{35},x_{36},...,x_{68})$ & $|Aut(G_i)   |$ & $\gamma$ & $\beta$ & Type \\ \hline
$1$ & $8$ &  $(0001001101110000000000101011001100)$ & $1$ & $\textbf{0}$ & $\textbf{218}$ & $W_{68,2}$  \\ \hline
$2$ & $7$ &  $(0110000010001000111000111000100010)$ & $1$ & $\textbf{1}$ & $\textbf{193}$ & $W_{68,2}$  \\ \hline
$3$ & $7$ &  $(1000100101011000011011110011000000)$ & $1$ & $\textbf{1}$ & $\textbf{195}$ & $W_{68,2}$  \\ \hline
$4$ & $7$ &  $(0101001010010010000100100101001001)$ & $1$ & $1$ & $233$ & $W_{68,2}$  \\ \hline
$5$ & $7$ &  $(0111010010001001001000000100101010)$ & $1$ & $\textbf{2}$ & $\textbf{193}$ & $W_{68,2}$  \\ \hline
$6$ & $7$ &  $(1100010011000010110111011101101111)$ & $1$ & $\textbf{2}$ & $\textbf{195}$ & $W_{68,2}$  \\ \hline
\end{tabular}}
\end{center}
\end{table}

\begin{table}[h!]\caption{New codes of length 68 from neighbours of $G_5$}
\begin{center}\scalebox{0.8}{
\begin{tabular}{ccccccc}
\hline
$H_{i}$ & $G_{i}$  & $(x_{35},x_{36},...,x_{68})$ & $|Aut(H_i)   |$ & $\gamma$ & $\beta$ & Type \\ \hline
$1$ & $5$ &  $(0010010110011000000010111001111110)$ & $1$ & $\textbf{1}$ & $\textbf{197}$ & $W_{68,2}$  \\ \hline
$2$ & $5$ &  $(0100001011001011101010110111011111)$ & $1$ & $\textbf{1}$ & $\textbf{199}$ & $W_{68,2}$  \\ \hline
$3$ & $5$ &  $(1101001011101101011111110111100111)$ & $1$ & $\textbf{2}$ & $\textbf{199}$ & $W_{68,2}$  \\ \hline
$4$ & $5$ &  $(0011000011001110011000001100000001)$ & $1$ & $\textbf{2}$ & $\textbf{191}$ & $W_{68,2}$  \\ \hline
$5$ & $5$ &  $(0001100100110010010101000111100100)$ & $1$ & $\textbf{2}$ & $\textbf{204}$ & $W_{68,2}$  \\ \hline
$6$ & $5$ &  $(1011101001000001101001010111011101)$ & $1$ & $\textbf{2}$ & $\textbf{218}$ & $W_{68,2}$  \\ \hline
\end{tabular}}
\end{center}
\end{table}

\begin{table}[h]\caption{Code of length 68 from the neighbours of $D_1$}
\begin{center}\scalebox{0.8}{
\begin{tabular}{ccccccc}
\hline
$I_{i}$ & $D_{i}$ & $(x_{35},x_{36},...,x_{68})$ & $|Aut(I_i)   |$ & $\gamma$ & $\beta$ & Type \\ \hline
$1$ & $1$ & $(1111000110110011110111001010111101)$ & $1$ & $5$ & $133$ & $W_{68,2}$  \\ \hline
\end{tabular}}
\end{center}
\end{table}

\begin{table}[h]\caption{Code of length 68 from the neighbours of $I_1$}
\begin{center}\scalebox{0.8}{
\begin{tabular}{ccccccc}
\hline
$J_{i}$ & $I_{i}$ & $(x_{35},x_{36},...,x_{68})$ & $|Aut(J_i)   |$ & $\gamma$ & $\beta$ & Type \\ \hline
$1$ & $1$ & $(0000100001011000111001010100001100$ & $1$ & $6$ & $141$ & $W_{68,2}$  \\ \hline
\end{tabular}}
\end{center}
\end{table}

\begin{table}[h]\caption{New codes of length 68 from the neighbours of $J_1$}\label{neighbours}
\begin{center}\scalebox{0.8}{
\begin{tabular}{ccccccc}
\hline
$K_{i}$ & $J_{i}$ & $(x_{35},x_{36},...,x_{68})$ & $|Aut(K_i)   |$ & $\gamma$ & $\beta$ & Type \\ \hline
$1$ & $1$ & $(1111111101001100010100001000010100)$ & $1$ & $\textbf{6}$ & $\textbf{131}$ & $W_{68,2}$  \\ \hline
$2$ & $1$ & $(0000001110010111101110011111001111)$ & $1$ & $\textbf{7}$ & $\textbf{158}$ & $W_{68,2}$  \\ \hline
\end{tabular}}
\end{center}

\bigskip

\caption{New codes of length 68 from the neighbours of $K_2$}\label{neighbours}
\begin{center}\scalebox{0.8}{
\begin{tabular}{ccccccc}
\hline
$L_{i}$ & $K_{i}$ & $(x_{35},x_{36},...,x_{68})$ & $|Aut(L_i)   |$ & $\gamma$ & $\beta$ & Type \\ \hline
$1$ & $2$ & $(0110111111010100011101010011010101)$ & $1$ & $\textbf{7}$ & $\textbf{155}$ & $W_{68,2}$  \\ \hline
$2$ & $2$ & $(0101010101010001001010011101110010)$ & $1$ & $\textbf{7}$ & $\textbf{156}$ & $W_{68,2}$  \\ \hline
$3$ & $2$ & $(0010011101010101010111011110110110)$ & $1$ & $\textbf{7}$ & $\textbf{157}$ & $W_{68,2}$  \\ \hline
$4$ & $2$ & $(1101111110110111001111110101101100)$ & $1$ & $\textbf{7}$ & $\textbf{159}$ & $W_{68,2}$  \\ \hline
$5$ & $2$ & $(1001011111000110001111101100101110)$ & $1$ & $\textbf{7}$ & $\textbf{160}$ & $W_{68,2}$  \\ \hline
$6$ & $2$ & $(1100000100100000010100101100011010)$ & $1$ & $\textbf{7}$ & $\textbf{162}$ & $W_{68,2}$  \\ \hline
$7$ & $2$ & $(1000010000010110000111110010011111)$ & $1$ & $\textbf{7}$ & $\textbf{164}$ & $W_{68,2}$  \\ \hline
$8$ & $2$ & $(0100001001101111111010000101010001)$ & $1$ & $\textbf{7}$ & $\textbf{165}$ & $W_{68,2}$  \\ \hline
$9$ & $2$ & $(0011101000100011011101001111101111)$ & $1$ & $\textbf{7}$ & $\textbf{167}$ & $W_{68,2}$  \\ \hline
\end{tabular}}
\end{center}
\end{table}

\newpage

\section{Conclusion}

In this work, we introduced a new construction that involved both block circulant matrices and a reverse circulant matrix. We demonstrated the relevance of this new construction by constructing many binary self-dual codes, including new self-dual codes of length $68$. To summarise the numerical results, we construct the following unknown $W_{68,1}$ code:
\[ \beta = \{317\}.\]
Furthermore, we construct the following unknown $W_{68,2}$ codes:
\begin{equation*}
\begin{split}
(\gamma =0,& \quad \beta =\{208,214,218\}), \\
(\gamma =1,& \quad \beta =\{179,191,193,195,197,199,202,210,211,229\}), \\
(\gamma =2,& \quad \beta =\{61,191,193,195,199,204,218\}), \\
(\gamma =6,& \quad \beta =\{131\}), \\
(\gamma =7,& \quad \beta =\{155,156,157,158,159,160,162,164,165,167\}) \\
\end{split}%
\end{equation*}

Regarding this construction, we were restricted to small group rings due to computational limitations. With a higher computational power,
it would be possible to investigate larger group rings which would yield more results. One could also consider other families of rings.

\end{document}